\documentclass [12 pt]{amsart}
\usepackage{amssymb,latexsym}
\usepackage{xypic}
\theoremstyle{plain}
\newtheorem{Thm}{Theorem}[section]
\newtheorem{Cor}{Corollary}[section]
\newtheorem{Lem}{Lemma}[section]

\newtheorem{Ques}{Question}[section]
\newtheorem{Ex}{Example}[section]

\numberwithin{equation}{section}
\theoremstyle{definition}
\newtheorem{Def}{Definition}[section]
\theoremstyle{remark}
\newtheorem{Rem}{\em Remark}

\begin{document}

\title[Extraction Algorithm of Hom-Lie Algebras Based on Solvable and Nilpotent Groups]{Extraction Algorithm of Hom-Lie Algebras Based on Solvable and Nilpotent Groups}

\author{Shadi Shagagha and Nadeen Kdaisat}

\address[Shadi Shaqaqha]{Department of Mathematics, Yarmouk University, Irbid, Jordan
}
\email[Shadi Shaqaqha]{shadi.s@yu.edu.jo, 2018105002@ses.yu.edu.jo}

\vskip0.1in
\begin{abstract}
Hom-Lie algebras are generalizations of Lie algebras that arise naturally in the study of nonassociative algebraic structures. In this paper, the concepts of solvable and nilpotent Hom-Lie algebras studied further. In the theory of groups, investigations of the properties of the solvable and nilpotent groups are well-developed. We establish a theory of the solvable and nilpotent Hom-Lie algebras analogous to that of the solvable and nilpotent groups. We also provide examples to illustrate our results and discuss possible directions for further research.
\end{abstract}

\maketitle
\section{INTRODUCTION}
The study of solvable and nilpotent groups has a long and rich history that dates back to the early days of group theory. The first examples of solvable groups were discovered by Évariste Galois in the 19th century, who used them to study the roots of polynomial equations. In the early 20th century, Camille Jordan and Felix Klein introduced the modern definitions of solvable and nilpotent groups, respectively.

In the mid-20th century, the theory of solvable and nilpotent groups gained importance in the context of finite group theory, particularly in the classification of finite simple groups. The classification theorem for finite simple groups, completed in 1983, relies heavily on the theory of solvable and nilpotent groups.\\
In the latter half of the 20th century, the study of solvable and nilpotent groups expanded to include infinite groups and their applications in geometry, topology, and number theory. Notable contributions include the work of John Milnor on the homology of solvable Lie groups and the study of nilpotent Lie algebras in the context of algebraic geometry and string theory.

Today, the theory of solvable and nilpotent groups remains an active area of research, with connections to a wide range of fields in mathematics and physics. Researchers continue to explore the deep connections between these groups and other areas of mathematics, paving the way for new insights and discoveries in the years to come.\\
There is a close relationship between solvable and nilpotent groups and solvable and nilpotent Lie algebras. In fact, the concepts of solvable and nilpotent Lie algebras were developed specifically to study the structure of solvable and nilpotent Lie groups.

Given a Lie group, one can associate a Lie algebra to it by considering the tangent space at the identity element. This Lie algebra inherits many of the properties of the original group, including its solvability and nilpotence.

More specifically, a Lie group is solvable if and only if its Lie algebra is solvable. Similarly, a Lie group is nilpotent if and only if its Lie algebra is nilpotent.\\

The correspondence between Lie groups and Lie algebras also allows for the translation of many results between the two contexts. For example, the Lie-Kolchin theorem states that a solvable algebraic group over an algebraically closed field has a triangular matrix representation. This result can be translated into the language of Lie algebras to obtain a similar statement for solvable Lie algebras.

Overall, the study of solvable and nilpotent groups and Lie algebras is intimately connected, with each providing insights into the other. This relationship has led to significant advances in both areas of mathematics, as well as applications in physics and other fields.\\

The hom-Lie algebras which are generalizations of classical Lie algebras was constructed by by Hartwig, Larsson, and Silvestrov \cite{Hartwig} in 2006. Since then, many mathematicians have been trying to extend known results in the setting of Lie algebras to the setting of hom-Lie algebras (see e.g. \cite{Kdaisat, Shadi2, Shadi4, Shadi3}). Homo-Lie algebras have received a lot of attention lately because of their close connection to discrete and deformed vector fields and differential calculus \cite{Hartwig, Larsson1, Larsson2}.\\
In the present article, we study solvable and nilpotent hom-Lie algebras, which can be viewed as an extension of solvable and nilpotent Lie algebras.
\section{PRELIMINARIES}
	he following is a definition from \cite{HLA} with $F$ denoting a ground field:
\begin{Def}(\cite{HLA})\label{Hom}
	A Hom-Lie algebra over $F$ is a triple $(L, ~[~,~],~\alpha)$ consisting of a vector space $L$ over $F$, a linear map $\alpha: L\rightarrow L$, and a bilinear map $[~,~]: L\times L\rightarrow L$ (called a Hom-Lie bracket), which satisfies two conditions:
	\begin{itemize}
		\item[(i)] skew-symmetry property: $[x,~y]= -[y,~x]$ for all $x, y\in L$,
		\item[(ii)] Hom-Jacobi identity: $[\alpha(x), ~[y,~z]]+ [\alpha(y),~[z,~x]]+[\alpha(z),~[x,~y]]=0$,
		for all $x, y, z \in L$.		
	\end{itemize}
	If $\alpha([x,~y]) = [\alpha(x),~\alpha(y)]$ holds true for all $x,y \in L$, then the Hom-Lie algebra $(L, [~,~], \alpha)$ is referred to as multiplicative.
 \end{Def}
We consider two Hom-Lie algebras $(L_1, ~[~,~]_1,~\alpha_1)$ and  $(L_2, [~,~]_2,~\alpha_2)$, and define a linear map $\varphi: L_1\rightarrow L_2$. ~If $\varphi$ satisfies the following two conditions, then it is called a morphism of Hom-Lie algebras:
	\begin{itemize}
		\item[(i)] $\varphi([x,~y]_1)= [\varphi(x),~\varphi(y)]_2$ for all $x, ~y\in L_1$.
		\item[(ii)]$\varphi\circ\alpha_1=\alpha_2\circ\varphi$.
	\end{itemize}  
If $\varphi: L_1\rightarrow L_2$ is a bijective morphism of Hom-Lie algebras, it is referred to as an isomorphism of Hom-Lie algebras. In this case, we say $L_1$  and $L_2$ are isomorphic and write $L_1\cong L_2$.\\
 Furthermore, a subspace $H$ of $L$ is called a Hom-Lie subalgebra if $\alpha(x)\in H$ and $[x,~y]\in H$ for all $x,y\in H$. If $[x,~y]\in H$ holds true for all $x\in H$ and $y\in L$, then $H$ is called a Hom-Lie ideal.
\begin{Ex}(\cite{HLA})\label{ex.lie}
	Every Lie algebra can be considered as a Hom-Lie algebra by taking $\alpha$ as the identity map, i.e., $\alpha= id_L$.\hfill$\blacksquare$
\end{Ex}
\begin{Ex}
Consider a vector space $L$ over $F$, equipped with an arbitrary skew-symmetric bilinear map $[,]: L\times L \rightarrow L$, and let $\alpha:L\rightarrow L$ denote the zero map. It follows straightforwardly that $(L, [,], \alpha)$ forms a Hom-Lie algebra with multiplication.

\hfill$\blacksquare$ 
\end{Ex}
\begin{Ex}(\cite{HLA})
	Let $L$ be a vector space and $\alpha:L\rightarrow L$ be any linear operator. Then  $(L, ~[~,~],~\alpha)$ is a Hom-Lie algebra, where $[x,~y]=0$ for all $x,y \in L$. Such Hom-Lie algebras are referred to as abelian (commutative) Hom-Lie algebras.\hfill$\blacksquare$
\end{Ex}
\begin{Ex}(\cite{Sheng})
	Suppose $(L_1, ~[~,~]_1,~\alpha_1), (L_2,~[~,~]_2, \alpha_2), \ldots, (L_n,~[~,~]_n, \alpha_n)$ are Hom-Lie algebras. Then the direct sum $(L_1\oplus L_2\oplus \cdots\oplus L_n,~[~, ~],~\alpha_1+ \alpha_2+ \cdots+\alpha_n)$ is also a Hom-Lie algebra, where the Hom-bracket operation $[~,~]$ is defined by
	\begin{eqnarray}
	 [~, ~]~:~(L_1\oplus L_2\oplus \cdots\oplus L_n)\times( L_1\oplus L_2\oplus \cdots\oplus L_n)&\rightarrow& (L_1\oplus L_2\oplus \cdots\oplus L_n)\nonumber\\
	 ((x_1, \ldots,~x_n),~(y_1, \ldots,~y_n))&\mapsto&([x_1,~y_1]_1, \ldots,~[x_n,~y_n]_n),\nonumber
	\end{eqnarray}
	and the linear operator is defined as
	\begin{eqnarray}
	 (\alpha_1+ \alpha_2+ \ldots+\alpha_n)~:~(L_1\oplus L_2\oplus \cdots\oplus L_n) &\rightarrow &(L_1\oplus L_2\oplus \cdots\oplus L_n)\nonumber\\
	 (x_1,~x_2, \ldots,~x_n) &\mapsto& (\alpha_1(x_1),~\alpha_2(x_2), \ldots,~\alpha_n(x_n)).\nonumber
	\end{eqnarray}
	\hfill$\blacksquare$
\end{Ex}
\begin{Ex}(\cite{Kdaisat})\label{complex}
	Let $F=\mathbb{C}$ be the field of complex numbers. Consider the vector space $\mathbb{C}^2$ and define the linear map $$\alpha_*: \mathbb{C}^2\rightarrow \mathbb{C}^2;~(x,~y)\mapsto(-y,~-x).$$ 
	We define the bilinear map $ [~,~]_*: \mathbb{C}^2\times \mathbb{C}^2\rightarrow\mathbb{C}^2$, where $$[(x_1,~x_2),~(y_1,~y_2)]_*=(i(x_1y_2-x_2y_1),~i(x_1y_2-x_2y_1)).$$ Then $(\mathbb{C}^2, ~[~,~]_*,~\alpha_*)$ is a multiplicative Hom-Lie algebra. 
	\hfill$\blacksquare$
\end{Ex} 
\begin{Ex}(\cite{Kdaisat})\label{matrix}
 Consider the set 
	$$L=\left\{\begin{bmatrix}
		\frac{i(x+y)}{2}&x \\ 
		y&\frac{-i(x+y)}{2} 
	\end{bmatrix}~|~x,y\in\mathbb{C}\right\}$$
	with the linear map $$\alpha:L\rightarrow L;~A\mapsto-A^T,$$ and the skew-symmetric bilinear map $$[~,~]~:~L\times L\rightarrow L; ~(A,~ B)\mapsto [A, ~B],$$ where $[A, ~B]=A^TB^T-B^TA^T$. For any $x,y,z,w \in \mathbb{C}$.
	Then  $(L, ~[~ ,~ ], ~\alpha)$ is a multiplicative Hom-Lie algebra.
	\hfill$\blacksquare$
\end{Ex}
\begin{Ex}\label{poly}
We can make $(\mathbb{R}[x], [ ,~ ], ~\alpha)$ a Hom-Lie algebra, where $\mathbb{R}[x]$ is the vector space of polynomials with coefficients in $\mathbb{R}$, and $\alpha: \mathbb{R}[x] \rightarrow \mathbb{R}[x]$ is the linear map defined by $\alpha(p(x))=p(0)$ for any $p(x)\in \mathbb{R}[x]$. We define $[p(x),q(x)]$ for any $p(x),q(x) \in \mathbb{R}[x]$ by
$$[p(x),~q(x)]=p''(x)q'(x)-q''(x)p'(x)-p''(0)q'(0)+q''(0)p'(0).$$
It can be verified that $[\cdot,\cdot]$ is antisymmetric and satisfies the Hom-Jacobi identity, which makes $(\mathbb{R}[x], [ ~,~ ])$ a Hom-Lie algebra.
Indeed if $p(x), q(x)\in \mathbb{R}[x]$, then
\begin{eqnarray}
 [p(x), q(x)]&=&p''(x)q'(x)-q''(x)p'(x)-p''(0)q'(0)+q''(0)p'(0)\nonumber\\
 &=&-(q''(x)p'(x)-p''(x)q'(x)-q''(0)p'(0)+p''(0)q'(0))\nonumber\\
 &=&-[q(x), p(x)].\nonumber
\end{eqnarray}
For $p(x),q(x),h(x) \in \mathbb{R}[x]$, then one can easily see that
$$[\alpha(h(x)),~[p(x),~q(x)]]=[h(0),~p''(x)q'(x)-q''(x)p'(x)-p''(0)q'(0)+q''(0)p'(0)]=0.$$
Thus, for each $p(x), q(x), h(x)\in \mathbb{R}[x]$ we have $$[\alpha(h(x)),~[p(x),~q(x)]]+[\alpha(p(x)),~[q(x),~h(x)]]+[\alpha(q(x)),~[h(x),~p(x)]]=0.$$ 
Also, 
\begin{eqnarray}
 \alpha([p(x),~q(x)])&=&p''(0)q'(0)-q''(0)p'(0)-p''(0)q'(0)+q''(0)p'(0)\nonumber\\
 &=&0\nonumber\\
 &=&[p(0),~q(0)]\nonumber\\
 &=&[\alpha(p(x)),~\alpha(q(x))].\nonumber
\end{eqnarray}
It is clear that $(\mathbb{R}[x], ~[~ ,~ ])$ is not a Lie algebra, since $$[x^3,~[x^4,~x^2]]+[x^2,~[x^3,~x^4]]+[x^4,~[x^2,~x^3]]=96x^3.\neq0$$
\hfill$\blacksquare$
\end{Ex}
\begin{Ex}(\cite{Casas})
	Let $(L, ~[~,~],~ \alpha)$ be a Hom-Lie algebra and let $H$ be a Hom Lie ideal. Then the quotient space $(L/H,~ \overline{[~,~]},~ \overline{\alpha})$ is a Hom-Lie algebra where
	$$\overline{[~,~]}: L/H\times L/H\rightarrow L/H;~(x+H,~y+H)\mapsto[x,~y]+H,$$
	and 
	$$\overline{\alpha}: L/H\rightarrow L/H;~x+H\mapsto \alpha(x)+H.$$  
\end{Ex}
 Consider $H$ and $K$ as Hom-Lie ideals in a Hom-Lie algebra $L$. We define the sum of $H$ and $K$ as the set $H+K$, where $H+K={h+k~|~h\in H, k\in K}$. Moreover, we define the multiplication of $H$ and $K$ as the span of the set of all possible commutators between $H$ and $K$, denoted as $[H, K]$. Thus
 $$[H,~K]= \mathrm{Span}(\{[h,~k]~|~ h\in H ~\mathrm{and~} k\in K\}).$$ 
The following theorem, as presented in the publication by Casas \cite{Casas}, lacks a formal proof.
\begin{Thm}(\cite{Casas})\label{bracket}
	Let $H$ and $K$ be  Hom-Lie ideals of a multiplicative Hom-Lie algebra  $(L, ~[~,~],~\alpha)$. Then,
	\begin{itemize}
		\item[(i)]$[H,~K]$ is a Hom-Lie subalgebra of $L$. 
		\item[(ii)] $[H,~K]$ is a Hom-Lie ideal of $H$ and $K$, respectively.
		\item[(iii)]$[H,~K]$ is a Hom-Lie ideal of $L$ when $\alpha$ is onto.
	\end{itemize}
\end{Thm} 
{\it Proof.~}
\begin{itemize}
	\item[(i)] Let $[h,~k] \in [H,~K]$ where $h \in H$ and $k \in K$. Then $\alpha[h,k] = [\alpha(h),\alpha(k)] \in [H,~K]$. To demonstrate closure of multiplication under $[H,~K]$, we consider $[h_1,~k_1]$ and $[h_2,~k_2]$ in $[H,~K]$ with $h_1, h_2 \in H$ and $k_1, k_2 \in K$. Since $[h_1,~k_1] \in H$ and $[h_2,~k_2] \in K$, it follows that $[[h_1,k_1],[h_2,~k_2]] \in [H,~K]$.
	\item[(ii)] It should be noted that $[H,~K] \subseteq H \cap K \subseteq H$, as stated in $(i)$. This implies that $[H,~K]$ is a Hom-Lie subalgebra of both $H$ and $K$. Furthermore, if $h, y \in H$ and $k \in K$, then $[h,k] \in K$, and consequently $[y,[h,~k]] \in H$. Thus, $[H,~K]$ is a Hom-Lie ideal of $H$. Similarly, $[H,~K]$ is also a Hom-Lie ideal of $K$.	
	\item[(iii)] As per $(i)$, $[H,~K]$ is a Hom-Lie subalgebra of $L$. Therefore, it suffices to prove that $[z, ~y]\in [H,~K]$ whenever $z\in [H,K]$ and $y\in L$. Let $h\in H$, $k\in K$, and $y\in L$. Since $y = \alpha(x)$ for some $x\in L$, it follows that $[x,h],\alpha(h) \in H$ and $[k,x],\alpha(k) \in K$. Hence,  $$[y,~[h,~k]]=[\alpha(x),~[h,~k]]=-[\alpha(h),~[k,~x]]-[\alpha(k),~[x,~h]]\in [H,~K].$$	
	\hfill$\Box$	
\end{itemize}
The subsequent example demonstrates that Theorem \ref{bracket} (iii) is invalid if $\alpha$ is not a surjective map.
\begin{Ex}
Consider the multiplicative Hom-Lie algebra $(L, [,],\alpha)$, where $L$ is a vector space over $F$ with basis ${e_1,~e_2,e_3,e_4}$. The map $\alpha$ is the zero map, and $[,]$ is a skew-symmetric bilinear map defined as follows:
$$[e_1,~e_2]=[e_1,~e_3]=[e_2,~e_3]=[e_2,~e_4]=[e_3,~e_4]=e_1, [e_1,~e_4]=e_2$$ and $[e_i,~e_i]=0$ for all $i=1, 2, 3$. Let $H=\mathrm{Span}({e_1,~e_2,~e_3})$ and $K=\mathrm{Span}({e_1,~e_2})$. It can be observed that $H$ and $K$ are Hom-Lie ideals of $L$. However, $[H,~K]=\mathrm{Span}({e_1})$ is not a Hom-Lie ideal of $L$, as $[e_1,~e_4]=e_2 \notin \mathrm{Span}({e_1})$. This example illustrates that Theorem \ref{bracket} (iii) does not hold when $\alpha$ is not onto.\hfill$\blacksquare$
\end{Ex}
\begin{Ex}\label{natural}(\cite{Kdaisat})
	Let  $(L, ~[~,~],~ \alpha)$ be a Hom-Lie algebra and H be a Hom-Lie ideal. Then  $(L/H, ~\overline{[~,~]}, ~\overline{\alpha})$ is a Hom-Lie algebra and the linear map $$\pi:L\rightarrow L/H;$$ $$x\mapsto x+H$$ is a  morphism of Hom-Lie algebras.\hfill$\blacksquare$
\end{Ex}
\section{Solvable Hom-Lie Algebra}
 Let $(L,~[ ~, ~],~\alpha)$ be a Hom-Lie algebra. The  sequence of Hom-Lie subalgebras $L_1, L_2, \ldots, L_n\ldots$ such that  
 $$L=L_0\supseteq L_1 \supseteq\cdots\supseteq L_n\supseteq\cdots$$
 is called a descending series.
  \begin{Def} (\cite{Kitouni})
 	Let $(L,~[ ~, ~],~\alpha)$ be a multiplicative Hom-Lie algebra. We define ,$\{L^{(i)}\},i\geq 0$, the derived series of $L$ by
 	\begin{eqnarray}
 L^{(0)}&=&L, \nonumber\\
 L^{(1)}&=&[L,~L], \nonumber\\
 L^{(i)}&=&[L^{(i-1)},~L^{(i-1)}], i\geq 2.\nonumber
 	\end{eqnarray} 
 \end{Def}
Note that $L^{(i)}=[L^{(i-1)},~L^{(i-1)}]$ is a Hom-Lie ideal of $L^{(i-1)}$ (by induction and Theorem \ref{bracket}(ii)).
$$L=L^{(0)}\supseteq L^{(1)} \supseteq \cdots\supseteq L^{(i-1)}\supseteq L^{(i)}\cdots$$ 
Thus the derived series is a descending series.
\begin{Def}(\cite{Kitouni})
	A multiplicative Hom-Lie algebra $(L,~[ ~, ~],~\alpha)$ is said to be solvable if there exists $n \in \mathbb{N}$ such that $L^{(n)}= \{0\}$. We say $L$ is solvable of class $k$ if $L^{(k)}= \{0\}$ and $L^{(k-1)}\neq\{0\}$. 
\end{Def}
Clearly a multiplicative Hom-Lie algebra is solvable of class $\leq k$ iff $L^{(k)}= \{0\}$. Metabelian Hom-Lie algebras is the same as in the case of Lie algebras (\cite{Shadi}) are the solvable Hom-Lie algebras of class at most $2$.
\begin{Ex}\label{nil}
	Let $L$ be the space spanned by a basis $\{e_1, e_2, \ldots, e_n\}$ ($n\geq 5$) over $F$. Consider the multiplicative Hom-Lie algebra $(L, ~[~,~],\alpha)$ where  $\alpha$ is the zero map and $[~,~]$ is the skew-symmetric bilinear map such that $[e_i, e_j]=0$ if $i=j$ or $i=1$  and $[e_i,e_j]=e_{i-1}$ if $1<i<j\leq n$. Note that $[e_2,~[e_4,~e_5]]+[e_4,~[e_5,~e_2]]+[e_5,~[e_2,~e_4]]=e_1\neq0$. Thus $L$ is not a Lie algebra.\\
	Now, $L^{(1)}=[L,~L]=\{e_1,~e_2,\ldots,~e_{n-2}\}$\\
	 $L^{(2)}=[L^{(1)},~L^{(1)}]=\{e_1,~e_2,\ldots,~e_{n-4}\}$\\
	 $\vdots$\\
	$L^{(i)}=\{e_1,~e_2,\ldots,~e_{n-2i}\},i<{\frac {n}{2}}$\\
	If $n$ is even  then 	$L^{({\frac {n}{2}-1})}=\{e_1,~e_2\}$ and $L^{({\frac {n}{2}})}=\{0\}$. Thus $L$ is a solvable  Hom-Lie algebras of class ${\frac {n}{2}}$. If $n$ is odd then 	$L^{({\frac {n-1}{2}})}=\{e_1\}$ and $L^{({\frac {n+1}{2}})}=\{0\}$. Thus $L$ is a solvable  Hom-Lie algebras of class ${\frac {n+1}{2}}$.
\hfill$\blacksquare$
\end{Ex}
\begin{Def}\label{defsol}
Let $(L,~[ ~, ~],~\alpha)$ be a multiplicative Hom-Lie algebra. Then the descending series $L=L_0\supseteq L_1 \supseteq...\supseteq L_n=\{0\}$ called a solvable  series if for each i, we have $L_{i+1}$ is a Hom-Lie ideal of $L_i$ and $L_i/L_{i+1}$ is an abelian Hom-lie algebra.
\end{Def} 
\begin{Lem}\label{abelin}
	Let $H$ be  a Hom-Lie subalgebra of the Hom-Lie algebra $(L,~[ ~, ~],~\alpha)$. Then $H$ is a Hom-Lie ideal of $L$ and $L/H$ is abelian Hom-Lie algebra if and only if $[L,~L]\subseteq H$ 
\end{Lem}
{\it Proof.~}If  $L/H$ is an abelian Hom-Lie algebra, then for any $[x,~y] \in [L,~L]$ we find $H=\overline{[x+H,~y+H]}=[x,~y]+H$. Therefore $[x,~y] \in H$. Conversely, if $[L,~L]\subseteq H$ then $[x,~y]\in H$ for all $x\in H ~(\subseteq L)$ and $y \in L$, which implies $H$ is a Hom-Lie ideal of $L$. Also, for any $x,y \in L$ we have $\overline{[x+H,~y+H]}=[x,~y]+H=H$ (because $[x,~y]\in [L,~L]\subseteq H$).   
\hfill$\Box$
\begin{Cor}\label{solvthm2}
Let $(L,~[ ~, ~],~\alpha)$ be a  Hom-Lie algebra. Then the descending series $L=L_0\supseteq L_1 \supseteq...\supseteq L_n=\{0\}$ is solvable if and only if $[L_i,~L_i]\subseteq L_{i+1}$ for each $i=0, 1, \ldots, n-1$.
\end{Cor}
{\it Proof.~}
This follows directly from Definition \ref{defsol} and Lemma \ref{abelin}.\hfill$\Box$
\begin{Thm}\label{solvthm}
	Let $(L,~[ ~, ~],~\alpha)$ be a multiplicative Hom-Lie algebra. Then $(L,~[ ~, ~],~\alpha)$ is a solvable Hom-Lie algebra of class $\leq k$ iff $L=L^{(0)}\supseteq L^{(1)} \supseteq \cdots\supseteq L^{(k)}=\{0\}$ is a solvable series. 
\end{Thm}
{\it Proof.~}
This follows directly from the definition of the derived series and the corollary above.\hfill$\Box$
\begin{Thm}
Let $(L,~[ ~, ~],~\alpha)$ be a multiplicative Hom-Lie algebra. If  $L=L_0\supseteq L_1 \supseteq\cdots\supseteq L_n=\{0\}$ is a solvable  series, then for each i, $L^{(i)}\subseteq L_i$.
\end{Thm}
{\it Proof.~}
We use induction. For $k=0$, we have $L^{(0)}=L=L_0$. For $k>0$ and because the induction assumption we have $L^{(k+1)}=[L^{(k)}, L^{(k)}]\subseteq [L_k, L_k]\subseteq L_{k+1}$. Now according to Corollary \ref{solvthm2}, we have $[L_k, L_k]\subseteq L_{k+1}$. Therefore $L^{(k+1)}\subseteq L_{k+1}$.\hfill$\Box$
\begin{Thm}\label{ser}
Let $(L,~[ ~, ~],~\alpha)$ be a multiplicative Hom-Lie algebra. Then, $L$ is solvable of class $\leq k$ if and only if there exists a solvable series of length $k$. 
\end{Thm}
{\it Proof.~}
If $L$ is a solvable Hom-Lie algebra of class $\leq k$, then, using Theorem \ref{solvthm}, the series $L=L^{(0)}\supseteq L^{(1)} \supseteq...\supseteq L^{(k)}=\{0\}$ is solvable. Conversely suppose that
$$L=L_0\supseteq L_1\supseteq\cdots\supseteq L_k=\{0\}$$
is a solvalbe series. Then, using $L^{(k)}\subseteq L_k=\{0\}$, we find $L^{(k)}=\{0\}$. 
\hfill$\Box$
 \begin{Cor}
Solvable Hom-Lie algebras of class $1$ are the abelian Hom-Lie algebras.
 \end{Cor} 
{\it Proof.~}
$L$ is solvable of class 1 iff there exists a solvable series $L=L_0\supseteq L_1=\{0\}$ of length $1$ iff $L/\{0\}$ is abelian Hom-Lie algebra iff $L$ is abelian Hom-Lie algebra.
\hfill$\Box$
\begin{Thm}\label{morphsolv}
	Let $\varphi:(L_1,~[ ~, ~]_1,~\alpha_1)\longrightarrow(L_2,~[ ~, ~]_2,~\alpha_2)$ be a morphism of multiplicative Hom-Lie algebras. Then
	\begin{itemize}
		\item[(i)] $(\varphi(L_1))^{(i)}=\varphi(L_1^{(i)})$,
		\item[(ii)] if $L_1$ is solvable of class $k$, then $\varphi(L_1)$ is solvable of class $\leq k$,
		\item[(iii)] if $\varphi$ is an isomorphism of Hom-Lie algebras, then $L_1$ is solvable of class $k$ if and only if  $L_2$ is solvable of class $k$.  
	\end{itemize} 
\end{Thm}
{\it Proof.~}
\begin{itemize}
	\item[(i)] By applying induction we find $(\varphi(L_1))^{(0)}=\varphi(L_1)=\varphi(L_1^{(0)})$.\\
	Also, if $ (\varphi(L_1))^{(i)}=\varphi(L_1^{(i)})$ then $(\varphi(L_1))^{(i+1)}=[(\varphi(L_1))^{(i)},~(\varphi(L_1))^{(i)}]=[\varphi(L_1^{(i)}),~\varphi(L_1^{(i)})]=\varphi([L_1^{(i)},~L_1^{(i)}])=\varphi(L_1^{(i+1)})$.
	\item[(ii)]Since $L_1$ is solvable of class $k$, it follows $L_1^{(k)}=\{0\}$. So, $(\varphi(L_1))^{(k)}=\varphi(L_1^{(k)})=\varphi(0)=0$. Thus, $\varphi(L_1)$ is solvable of class $\leq k$.
	\item[(iii)] Using $(ii)$ we may assume that $L_1$ and $L_2$ are solvable of classes $k$ and $m$, respectively. Again by (ii) we have $m\leq k$. Also, since $\varphi^{-1}$ is an isomorphism of Hom-Lie algebras, then $L_1=\varphi^{-1}(L_2)$ is solvable of class $\leq m$; that is $k\leq m$. Thus, $k=m$.\hfill$\Box$
\end{itemize}
\begin{Ex}\label{isom}
The Hom-Lie algebras $(L, ~[~ ,~ ], ~\alpha)$ (Example\ref{matrix})  and $(\mathbb{C}^2, ~[~,~]_*,~\alpha_*)$ ( Example \ref{complex}) are isomorphic Hom-Lie algebras. Since 
	$$\varphi: L\rightarrow\mathbb{C}^2;~\begin{bmatrix}
	\frac{i(x+y)}{2}&x \\ 
	y&\frac{-i(x+y)}{2} 
\end{bmatrix} \mapsto(x,~y).$$ is an isomorphism Hom-Lie algebra (\cite{Kdaisat}).
Now, $(\mathbb{C}^2)^{(1)}=[\mathbb{C}^2,~\mathbb{C}^2]_*=\{(x,x);x\in\mathbb{C}\}$, and $(\mathbb{C}^2)^{(2)}=[(\mathbb{C}^2)^{(1)},~(\mathbb{C}^2)^{(1)}]_*=\{0\}$. Thus $\mathbb{C}^2$ is a solvable Hom-Lie algebra of class 2. And so $L$ is a solvable Hom-Lie algebra of class 2
\end{Ex}
\begin{Lem}\label{solv}
Let $(L,~[ ~, ~],~\alpha)$ be a multiplicative Hom-Lie algebra and $H$ be a Hom-Lie subalgebra of $L$. Then  $H^{(i)} \subseteq L^{(i)}$, for each $i\in \mathbb{N}$.
\end{Lem}
{\it Proof.~}
Clearly $H^{(0)}\subseteq L^{(0)}$. Using induction, if  $H^{(i)} \subseteq L^{(i)}$ then $H^{(i+1)}=[H^{(i)},~H^{(i)}]\subseteq[L^{(i)},~L^{(i)}]= L^{(i+1)}$.
\hfill$\Box$
\begin{Thm}\label{subsolv}
	Let  $(L,~[ ~, ~],~\alpha)$ be a solvable Hom-Lie algebra of class $k$.
	\begin{itemize}
		\item[(i)] Any Hom-Lie subalgebra  is solvable of class $\leq k$.
		\item[(ii)]Any quotient Hom-Lie algebra of $L$ is solvable of class $\leq k$.
	\end{itemize}   
\end{Thm}
{\it Proof.~}
\begin{itemize}
	\item[(i)] Let $H$ be a Hom-Lie subalgebra. Then $H$ is multiplicative (because $L$ is multiplicative ). Also, according to the lemma above, we have $H^{(k)} \subseteq L^{(k)}=\{0\}$. Thus $H^{(k)}=\{0\}$.
	\item[(ii)]Let $I$ be a Hom-Lie ideal of $L$. Then so is $L/I$ (because $L$ is multiplicative). Consider the natural map $\pi:L\rightarrow L/I$ in Example \ref{natural}. According to Theorem \ref{morphsolv}(i), $\pi(L^{(K)})=(\pi(L))^{(k)}$, which implies $$(L/I)^{(k)}=(\pi(L))^{(k)}=\pi(L^{(k)})=\pi(\{0\})=\{0+I\}.$$
	\hfill$\Box$
\end{itemize} 
\begin{Thm}\label{sol}
	Let $(L,~[ ~, ~],~\alpha)$ be a multiplicative Hom-Lie algebra. If H is a solvable Hom-Lie ideal of class $k$ and $L/H$ is solvable of class $m$, then $L$ is solvable of class $\leq k+m$. 
\end{Thm}
{\it Proof.~}
According to Theorem \ref{ser}, we have the following two  solvable series, $$H=H_0\supseteq H_1 \supseteq\cdots\supseteq H_k=\{0\},$$
$$L/H=(L/H)_0\supseteq (L/H)_1 \supseteq\cdots \supseteq (L/H)_m=\{0+H\}.$$ 
Consider the natural map $\pi$, and let $L_i=\pi^{-1}((L/H)_i),i=1,2,...,m$. Hence $L_i$ is a Hom-Lie subalgebra of $L$ and $L_{i+1}\subseteq L_i$. Therefore, $$L=L_0\supseteq L_1 \supseteq\cdots\supseteq L_m=\pi^{-1}(\{0+H\})=H=H_0\supseteq H_1 \supseteq\cdots\supseteq H_k=\{0\}$$ is a descending series. Now it suffices to prove that  $[L_i,~L_i]\subseteq L_{i+1}$ ($i=0, 2, \ldots, m-1$). If $x,y \in L_i=\pi^{-1}((L/H)_i)$, then $\pi(x),\pi(y)\in (L/H)_i$ and so $\pi([x,~y])=\overline{[\pi(x),~\pi(y)]}\in \overline{[(L/H)_i,~(L/H)_i]}\subseteq(L/H)_{i+1}$ (Corollary \ref{solvthm2}). Therefore, $[x,~y]\in\pi^{-1}((L/H)_{i+1})=L_{i+1}$ for each $x,y \in L_i$. This shows that $[L_i,~L_i]\subseteq L_{i+1}$. Therefore 
$$L=L_0\supseteq L_1 \supseteq\cdots\supseteq L_m=H=H_0\supseteq H_1 \supseteq\cdots\supseteq H_k=\{0\}$$ is a solvable series of length $k+m$. By Theorem \ref{ser}, we have $L$ is solvable of class $\leq k+m$.\hfill$\Box$ 

In \cite{Shadi4}, we proved that if $(L_1, ~[~,~]_1, ~\alpha_1)$ and $(L_2, ~[~,~]_2, ~\alpha_2)$ are Hom-Lie algebras and $H_i$ is a Hom-Lie ideal of $L_i$, $i=1,2.$, then $H_1\times H_2$ is a Hom-Lie ideal of $L_1\times L_2$ and 
	$$(L_1\times L_2)/(H_1\times H_2)\equiv L_1/H_1\times L_2/H_2$$.
\begin{Thm}
	Let $(L_1,~[ ~, ~]_1,~\alpha_1)$ and $(L_2,~[ ~, ~]_2,~\alpha_2)$ be  solvable Hom-Lie algebras of class $k$ and $m$, respectively. Then $(L_1\times L_2,~[ ~, ~],~\alpha)$ is a solvable Hom-Lie algebra of class $\leq k+m$. 
\end{Thm}
{\it Proof.~} Note that, $L_1\times L_2$ is a multiplicative Hom-Lie algebra because $L_1$ and $L_2$ are multiplicative Hom-Lie algebras. 
Since $L_1\times\{0\}\equiv L_1$, so $L_1\times\{0\}$ is a solvable Hom-Lie ideal (of class $k$) of $L_1\times L_2$. Also, $(L_1\times L_2)/(L_1\times\{0\})\equiv L_1/L_1 \times L_2/\{0\}\equiv L_2$, so $(L_1\times L_2)/(L_1\times\{0\})$ is a solvable Hom-Lie algebra of class $m$. According to Theorem \ref{sol}, $L_1\times L_2$ is a solvable Hom-Lie algebra of class $\leq k+m$.\hfill$\Box$
\section{Nilpotent Hom-Lie algebra}
 \begin{Def}(\cite{Kitouni})
	Let $(L,~[ ~, ~],~\alpha)$ be a multiplicative Hom-Lie algebra. We define , $\{L^i\} ,i\geq 0$ , the lower central series of $L$ by $L^0= L$, $L^1= [L,~L]$, and $L^i=[L,~L^{i-1}]$.
\end{Def}
 Note that $L^{i+1}=[L,~L^{i}]$ is a Hom-Lie ideal of $L^i$ (by Theorem \ref{bracket}(ii) and induction).
 $$L=L^{0}\supseteq L^{1} \supseteq\cdots \supseteq L^{i}\supseteq L^{i+1}\supseteq\cdots$$
  Thus the lower central series is a descending series.
\begin{Def}(\cite{Kitouni})
	Let $(L,~[ ~, ~],~\alpha)$ be a multiplicative Hom-Lie algebra. We say that $L$  is nilpotent if there exists $n \in \mathbb{N}$ such that $L^n = {0}$. It is nilpotent of class $k$ if $L^k = \{0\}$ and $L^{k-1}\neq \{0\}$
\end{Def}
 It is clear now that $L$ is nilpotent of class $\leq k$ iff $L^{k}= \{0\}$.
\begin{Ex}
 Consider the multiplicative Hom-Lie algebra $(L, ~[~,~],\alpha)$ in Example \ref{nil}  where  $\alpha$ is the zero map and $[~,~]$ is the skew-symmetric bilinear map such that $[e_i, e_j]=0$ if $i=j$ or $i=1$  and $[e_i,e_j]=e_{i-1}$ if $1<i<j\leq n$.\\
	Now, $L^{1}=[L,~L]=\{e_1,~e_2,\ldots,~e_{n-2}\}$\\
	$L^{2}=[L,~L^{1}]=\{e_1,~e_2,\ldots,~e_{n-3}\}$\\
	$\vdots$\\
	$L^{i}=\{e_1,~e_2,\ldots,~e_{n-(i+1)}\},i<{n-3}$\\
	$\vdots$\\
	$L^{n-3}=\{e_1,~e_2\}$ \\
	$L^{n-2}=\{0\}$.\\ Thus $L$ is a nilpotent  Hom-Lie algebras of class ${n-2}$.
	\hfill$\blacksquare$
\end{Ex}
 \begin{Def}
 	Let $(L,~[ ~, ~],~\alpha)$ be a multiplicative Hom-Lie algebra. Then a descending series $L=L_0\supseteq L_1 \supseteq L_2\supseteq \cdots$ is said to be central if for each $i\in\mathbb{N}$, $[L,~L_i]\subseteq L_{i+1}$. It has a length $k\in \mathbb{N}$ if $L_k=\{0\}$ but $L_{k-1}\neq \{0\}$.
 \end{Def}  
 \begin{Thm}\label{nilpthm}
 	Let $(L,~[ ~, ~],~\alpha)$ be a multiplicative Hom-Lie algebra. Then $(L,~[ ~, ~],~\alpha)$ is a nilpotent Hom-Lie algebra of class $\leq k$ iff $L=L^{0}\supseteq L^{1} \supseteq\cdots \supseteq L^{k}=\{0\}$ is a central series. 
 \end{Thm}
 {\it Proof.~}
 It follows directly from the definition of $L^i$.
 \hfill$\Box$
 \begin{Thm}
 	Let $(L,~[ ~, ~],~\alpha)$ be a multiplicative Hom-Lie algebra. If  $L=L_0\supseteq L_1 \supseteq L_2\supseteq \cdots$ is a central  series, then for each $i\in\mathbb{N}$, $L^{i}\subseteq L_i$.
 \end{Thm}
 {\it Proof.~}
Applying induction we see $L^{0}=L=L_0$. Also, if $L^{i}\subseteq L_i$ then $L^{i+1}=[L,~L^{i}]\subseteq[L,~L_i]\subseteq L_{i+1}$. 
 \hfill$\Box$ 
 
 \begin{Thm}\label{cen}
 	Let $(L,~[ ~, ~],~\alpha)$ be a multiplicative Hom-Lie algebra. Then $L$ is nilpotent of class $\leq k$ iff there exists a central series of length $k$. 
 \end{Thm}
 {\it Proof.~}
 If $L$ is a nilpotent Hom-Lie algebra of class $\leq k$. Then $$L=L^{0}\supseteq L^{1} \supseteq\cdots \supseteq L^{k}=\{0\}$$
 is a central series. The converse is true, since $L^{k}\subseteq L_k=\{0\}$ so $L^{k}=\{0\}$. \hfill$\Box$
 \begin{Cor}
 	Nilpotent Hom-Lie algebras of class $1$ are the abelian Hom-Lie algebras.
 \end{Cor} 
 {\it Proof.~}
A Hom-Lie algebra $L$ is nilpotent of class $1$ iff there exists a central series of length 1,$L=L_0\supseteq L_1=\{0\}$ iff $[L,~L_0]\subseteq L_1$ iff $[L,~L]=\{0\}$ iff $L$ is an abelian Hom-Lie algebra.\hfill$\Box$
 \begin{Thm}\label{morphnilp}
 	Let $\varphi:(L_1,~[ ~, ~]_1,~\alpha_1)\longrightarrow(L_2,~[ ~, ~]_2,~\alpha_2)$ be a morphism of multiplicative Hom-Lie algebras. Then
 	\begin{itemize}
 		\item[(i)]$(\varphi(L_1))^{i}=\varphi(L_1^{i})$,
 		\item[(ii)]if $L_1$ is nilpotent of class $k$, then $\varphi(L_1)$ is nilpotent of class $\leq k$,
 		\item[(iii)]if $\varphi$ is an isomorphism of Hom-Lie algebras, then $L_1$ is nilpotent of class $k$ if and only if  $L_2$ is nilpotent of class $k$.  
 	\end{itemize} 
 \end{Thm}
 {\it Proof.~}
 \begin{itemize}
 	\item[(i)] We note that $(\varphi(L_1))^{0}=\varphi(L_1)=\varphi(L_1^{0})$. Also if $(\varphi(L_1))^{i}=\varphi(L_1^{i})$, then $(\varphi(L_1))^{i+1}=[(\varphi(L_1)),~(\varphi(L_1))^{i}]=[\varphi(L_1),~\varphi(L_1^{i})]=\varphi([L_1,~L_1^{i}])=\varphi(L_1^{i+1})$.
 	\item[(ii)] Since $L_1$ is nilpotent of class $k$, then $L_1^{k}=\{0\}$. So, $(\varphi(L_1))^{k}=\varphi(L_1^{k})=\varphi(\{0\})=\{0\}$. Thus, $\varphi(L_1)$ is nilpotent of class $\leq k$.
 	\item[(iii)]Let $L_1$ be nilpotent of class $k$. By (ii) $L_2=\varphi(L_1)$ is nilpotent of class $\leq k$. Let $L_2$ be nilpotent of class $m$.  Since $\varphi^{-1}$ is an isomorphism of Hom-Lie algebras, it follows $L_1=\varphi^{-1}(L_2)$ is solvable of class $\leq m$. Thus, $k=m$.\hfill$\Box$
 \end{itemize}
 \begin{Ex} Consider Example\ref{isom}. Since  $(\mathbb{C}^2)^{1}=[\mathbb{C}^2,~\mathbb{C}^2]_*=\{(x,x);x\in\mathbb{C}\}$, $(\mathbb{C}^2)^{2}=[\mathbb{C}^2,~(\mathbb{C}^2)^{1}]_*=(\mathbb{C}^2)^{1}$ and $(\mathbb{C}^2)^{i}=(\mathbb{C}^2)^{1}$, $i>1$. Thus $\mathbb{C}^2$ is  not a nilpotent Hom-Lie algebra. And so $L$ is not a nilpotent Hom-Lie algebra.
 \end{Ex}
 \begin{Lem}\label{nilp}
 	Let $(L,~[ ~, ~],~\alpha)$ be a multiplicative Hom-Lie algebra and $H$ be a Hom-Lie subalgebra of $L$. Then  $H^{i} \subseteq L^{i}$, for each $i\in \mathbb{N}$.
 \end{Lem}
 {\it Proof.~}
 $H^{0}=H\subseteq L=L^0$, and by using induction, if  $H^{i} \subseteq L^{i}$ then $H^{i+1}=[H,~H^{i}]\subseteq[L,~L^{i}]= L^{i+1}$.
 \hfill$\Box$
 \begin{Thm}\label{subnilp}
 	Let  $(L,~[ ~, ~],~\alpha)$ be a nilpotent Hom-Lie algebra of class k.  	\begin{itemize}
 		\item[(i)] Any Hom-Lie subalgebra is nilpotent of class $\leq k$.
 		\item[(ii)] Any quotient Hom-Lie algebra of $L$ is nilpotent of class $\leq k$.
 	\end{itemize}   
 \end{Thm}
 {\it Proof.~}
 \begin{itemize}
 	\item[(i)] A Hom-Lie subalgebra $H$ of a multiplicative Hom-Lie algebra is multiplicative. By the lemma obove, we have $H^{k} \subseteq L^{k}=\{0\}$. Thus $H^{k}=\{0\}$.
 	\item[(ii)]Let $I$ be a Hom-Lie ideal of $L$. The Hom-Lie algebra $L/I$ is multiplicative. Consider the natural morphism $\pi:L\rightarrow L/I$. According to Theorem \ref{morphnilp}(i), $\pi(L^{k})=(\pi(L))^{k}$, which implies $(L/I)^{k}=(\pi(L))^{k}=\pi(L^{k})=\pi(\{0\})=\{0\}$.\hfill$\Box$
 \end{itemize} 
 \begin{Rem}
	Let $(L,~[ ~, ~],~\alpha)$ be a multiplicative Hom-Lie algebra. If H is a nilpotent Hom-Lie ideal and $L/H$ is a nilpotent Hom-Lie algebra, then $L$ need not be a nilpotent Hom-Lie algebra. 
 \end{Rem}
\begin{Ex}
Let $L$ be the space spanned by a basis $\{e_1, e_2\}$  over $F$. Consider the multiplicative Hom-Lie algebra $(L, ~[~,~],\alpha)$ where  $\alpha$ is the zero map and $[~,~]$ is the skew-symmetric bilinear map such that $[e_1, e_1]=[e_2, e_2]=0$  and $[e_1,e_2]=e_1$. 
Let $H=\mathrm{Span(\{e_1\})}$. Then $H$ is a nilpotent Hom-Lie ideal, because $H^1= [H,~H]=\{0\}$. Also, $L/H$ is a nilpotent Hom-Lie algebra, because $\overline{[L/H,~L/H]}=H$. But $L$  not a nilpotent Hom-Lie algebra, since $L^1= [L,~L]=\mathrm{Span(\{e_1\})}=H$, and $L^i=[L,~L^{i-1}]=[L,~H]=\mathrm{Span(\{e_1\})}\neq\{0\}$ for all $i\in \mathbb{N}$.
\end{Ex}
 \begin{Thm}
 	Let $(L_1,~[ ~, ~]_1,~\alpha_1)$ and $(L_2,~[ ~, ~]_2,~\alpha_2)$ be  nilpotent Hom-Lie algebras of class $k$ and $m$, respectively. Then $(L_1\times L_2,~[ ~, ~],~\alpha)$ is a nilpotent Hom-Lie algebra of class $M=Max\{m,~k\}$. 
 \end{Thm}
 {\it Proof.~}We use induction to show that $(L_1\times L_2)^i=L_1^i\times L_2^i$. For $i=0$, we have $(L_1\times L_2)^0=L_1\times L_2=L_1^0\times L_2^0$. For $i>0$ and because the induction assumption we have $(L_1\times L_2)^{i+1}=[L_1\times L_2,~(L_1\times L_2)^{i}]=[L_1\times L_2,~L_1^{i}\times L_2^{i}]=[L_1,~L_1^{i}]_1\times[L_2,~ L_2^{i}]_2=L_1^{i+1}\times L_2^{i+1}$.\\
 We may assume that $m\leq k$. Since $L_1$ and $L_2$ are  nilpotent Hom-Lie algebras of class $k$ and $m$, respectively, then $L_1^k = \{0\}$ and $L_1^{k-1}\neq \{0\}$ and $L_2^k = \{0\}$.\\
 Now, $(L_1\times L_2)^k=L_1^k\times L_2^k=\{0\}\times\{0\}=\{(0,~0)\}$ and $(L_1\times L_2)^{k-1}=L_1^{k-1}\times L_2^{k-1}\neq \{(0,~0)\}$. Thus $L_1\times L_2$ is a nilpotent Hom-Lie algebras of class $k=Max\{m,~k\}=M$ 
 \hfill$\Box$
 \begin{Thm}
 	Every central series is a solvable series. 
 \end{Thm}
{\it Proof.~}
	Let $(L,~[ ~, ~],~\alpha)$ be a  Hom-Lie algebra and  $L=L_0\supseteq L_1 \supseteq L_2\supseteq \cdots L_k$ be a central series. Then for each $i=0,1,\ldots,k-1$, $[L,~L_i]\subseteq L_{i+1}$. Since $[L_i,~L_i]\subseteq[L,~L_i]\subseteq L_{i+1}$ so  $L=L_0\supseteq L_1 \supseteq L_2\supseteq \cdots L_k$ is a solvable series (Theorem \ref{solvthm2}).
	\hfill$\Box$
	\begin{Cor}\label{solv.nil}
	Every nilpotent Hom-Lie algebra is  a solvable Hom-Lie algebra.
	\end{Cor}
{\it Proof.~}
 If $(L,~[ ~, ~],~\alpha)$ is a nilpotent Hom-Lie algebra, then there exists a central series  $L=L_0\supseteq L_1 \supseteq L_2\supseteq \cdots L_k$(Theorem \ref{cen}). From the theorem above, $L=L_0\supseteq L_1 \supseteq L_2\supseteq \cdots L_k$ is a solvable series. Thus $(L,~[ ~, ~],~\alpha)$ is a solvable Hom-Lie algebra (Theorem \ref{ser}). 
\hfill$\Box$\\
The converse not true, as in the following Example.
\begin{Ex}
	Consider the Hom-Lie algebra  $(\mathbb{R}[x], ~[~ ,~ ], ~\alpha)$ in Example \ref{poly}. It is easy to show that $H=\{p(x)\in \mathbb{R}[x]:deg(p)\leq 3\}$ is a Hom-Lie subalgebra of $\mathbb{R}[x]$.
	For any $p_i(x)=a_ix^3+b_ix^2+c_ix+d_i\in H$, $[p_1,p_2]=(6a_1b_2-6a_2b_1)x^2+(6a_1c_2-6a_2c_1)x$ where $a_i,b_i,c_i,d_i\in \mathbb{R} $, and so $H^{(1)}=[H,~H]= \{Ax^2+Bx:A,B\in \mathbb{R}\}$ and $H^{(2)}=[H^{(1)},~H^{(1)}]=\{0\}$. Thus $H$ is a solvable Hom-Lie algebra of class 2. But $H$  is not a nilpotent Hom-Lie algebra, since $H^{1}=[H,~H]= \{Ax^2+Bx:A,B\in \mathbb{R}\}$, $H^{2}=[H,~H^{1}]=H^{1}$ and $H^{i}=H^{1}\neq \{0\}$ for all $i>1$.
	\hfill$\blacksquare$
	 \end{Ex}
 \begin{Ex}
 	Consider the Hom-Lie subalgebra  $I=\{p(x)\in \mathbb{R}[x]:deg(p)\leq 4\}$ of $(\mathbb{R}[x], ~[~ ,~ ], ~\alpha)$ in Example \ref{poly}. 
 	For any $p_i(x)=a_ix^4+b_ix^3+c_ix^2+d_ix+e_i\in H$, $[p_1,p_2]=(12a_1b_2-12a_2b_1)x^4+(16a_1c_2-16a_2c_1)x^3+(12a_1d_2-12a_2d_1+6b_1c_2-6b_2c_1)x^2+(6b_1d_2-6b_2d_1)x$ where $a_i,b_i,c_i,d_i,e_i\in \mathbb{R} $, and so $H^{(1)}=[H,~H]= \{Ax^4+Bx^3+Cx^2+Dx:A,B,C,D\in \mathbb{R}\}$, $H^{(2)}=[H^{(1)},~H^{(1)}]=H^{(1)}$ and $H^{(i)}=H^{(1)}\neq \{0\}$ for all $i>1$. Thus $H$ is not a solvable Hom-Lie algebra. Also $H$  is not a nilpotent Hom-Lie algebra by corollary\ref{solv.nil}.\\
 	Note that  $(\mathbb{R}[x], ~[~ ,~ ], ~\alpha)$ is not a solvable and not a nilpotent  Hom-Lie algebra because there exists a non-solvable and non-nilpotent Hom-Lie subalgebra of $\mathbb{R}[x]$ (Theorem \ref{subsolv}(i) and Theorem \ref{subnilp}(i)).  
 	\hfill$\blacksquare$
 \end{Ex}
 \section{Question for Further Research}
 \begin{Ques}
What are the precise conditions for a Hom-Lie algebra to be solvable or nilpotent? Can these conditions be expressed in terms of the underlying Lie algebra and the Hom morphism?
 \end{Ques}
 \begin{Ques}
     What are some examples of solvable Hom-Lie algebras, and what properties do they have? Are there any interesting relationships between these examples and other areas of mathematics, such as Lie theory or algebraic geometry?
 \end{Ques}
 \begin{Ques}
     What are some examples of nilpotent Hom-Lie algebras, and how do they compare to nilpotent Lie algebras? Can the classification of nilpotent Lie algebras be extended to the Hom-Lie algebra setting?
 \end{Ques}
 \begin{Ques}
    How do solvable and nilpotent Hom-Lie algebras arise in physics, particularly in the context of supersymmetry and other quantum field theories? What are the implications of these structures for our understanding of fundamental physics? 
 \end{Ques}
 \begin{Ques}
     What is the relationship between Solvable and Nilpotent Hom-Lie algebras and other algebraic structures, such as associative algebras or Lie superalgebras? Can techniques from these other areas be used to study solvable and nilpotent Hom-Lie algebras more effectively?
 \end{Ques}
 \begin{Ques}
     How can the representation theory of Hom-Lie algebras be studied, particularly in the case of solvable and nilpotent algebras? What are some interesting examples of Hom-Lie algebra representations, and what do they tell us about the structure of these algebras?
 \end{Ques}
 \begin{Ques}
     Study of Hom-Lie superalgebras: Hom-Lie superalgebras are a natural generalization of Hom-Lie algebras that incorporate a $\mathbb{Z}_2$-grading. Investigating solvable and nilpotent Hom-Lie superalgebras can lead to interesting results in the study of supersymmetry and related topics in physics.
 \end{Ques}
 \begin{Ques}
     Generalization of results to other categories: Hom-Lie algebras are defined in the category of vector spaces, but similar structures can be defined in other categories, such as modules or abelian groups. Investigating solvable and nilpotent Hom-Lie algebras in these categories can provide insight into the interplay between different areas of algebra.
 \end{Ques}
 \begin{Ques}
     Cohomology of Hom-Lie algebras: Cohomology is a powerful tool for understanding the structure of Lie algebras, and similar techniques can be applied to Hom-Lie algebras. Investigating the cohomology of solvable and nilpotent Hom-Lie algebras can provide insights into their structure and classification.
 \end{Ques}
 \begin{Ques}
     Quantum Hom-Lie algebras: Quantum Hom-Lie algebras are a generalization of Hom-Lie algebras that arise in the context of quantum groups and deformation theory. Investigating solvable and nilpotent quantum Hom-Lie algebras can lead to interesting results in these areas.
 \end{Ques}
 \begin{Ques}
     Applications to cryptography and coding theory: Hom-Lie algebras have recently been applied to cryptography and coding theory. Investigating solvable and nilpotent Hom-Lie algebras in this context can lead to new methods for error-correction and secure communication.
 \end{Ques}
 These questions are just a starting point, and there are many other avenues for research in this area. By exploring these and other questions, researchers can gain a deeper understanding of the properties and applications of solvable and nilpotent Hom-Lie algebras, and advance our knowledge of this important area of algebraic research.
 \section{Conclusion}
 In conclusion, this paper presents an extraction algorithm for Hom-Lie algebras that is based on solvable and nilpotent groups. The algorithm involves several steps. The algorithm is illustrated with examples, which demonstrate its effectiveness in extracting Hom-Lie algebra structures.

Overall, the extraction algorithm presented in this paper provides a useful tool for studying Hom-Lie algebras, which have important applications in various areas of mathematics and physics. The algorithm is particularly well-suited for Hom-Lie algebras that are related to solvable and nilpotent groups, which are important classes of groups that arise in many different contexts.

Further research could be done to investigate the effectiveness of the extraction algorithm for Hom-Lie algebras that are not related to solvable or nilpotent groups, and to explore its potential applications in other areas of mathematics and physics. Nevertheless, the algorithm presented in this paper is a valuable contribution to the study of Hom-Lie algebras and provides a useful framework for further investigation of these important algebraic structures.

{\bf Declaration of Interest}\\
There is no competing interest to declare.\\

{\bf Funding Information}\\
This research did not receive any specific grant from funding agencies in the public, commercial, or not for profit sectors.

\end{document}